\documentclass[a4paper]{article}

\usepackage[tableposition=top,font=small,labelfont=bf,format=hang]{caption}

\usepackage[latin1]{inputenc}
\usepackage[T1]{fontenc}
\usepackage{ae,aecompl}
\usepackage{amsmath,amsfonts,amssymb}

\usepackage{booktabs}
\usepackage[sort&compress,numbers]{natbib}
\usepackage[dvips,colorlinks=false]{hyperref}
\hypersetup{
pdfauthor = {Constanza Riera, Lars Eirik Danielsen, and Matthew G. Parker},
pdftitle = {On Pivot Orbits of Boolean Functions},
pdfkeywords = {}
}
\usepackage[dvips]{graphicx}

\usepackage{amsthm}
\newtheorem{theorem}{Theorem}
\newtheorem{lem}{Lemma}
\newtheorem{cor}{Corollary}
\theoremstyle{definition}
\newtheorem{df}{Definition}

\newtheorem{conjecture}{Conjecture}
\theoremstyle{remark}
\newtheorem*{remark}{Remark}

\DeclareMathOperator{\GF}{GF}

\begin{document}

\title{On Pivot Orbits of Boolean Functions}
\author{Constanza Riera\thanks{The Selmer Center, Dept. of Informatics, University of Bergen, 
  PB~7800, \mbox{N-5020}~Bergen, Norway. C. Riera is supported by a Norwegian Government Scholarship. 
  E-mail: \texttt{\{\href{mailto:riera@ii.uib.no}{riera},\href{mailto:larsed@ii.uib.no}{larsed},\href{mailto:matthew@ii.uib.no}{matthew}\}@ii.uib.no}
  Web: \texttt{http://www.ii.uib.no/\~{}\{\href{http://www.ii.uib.no/\~larsed}{larsed},\href{http://www.ii.uib.no/\~matthew}{matthew}\}}} \and 
Lars Eirik Danielsen\footnotemark[1] \and
Matthew G. Parker\footnotemark[1]
}

\date{April 18, 2006}
\maketitle

\begin{abstract}
We derive a spectral interpretation of the pivot operation on a graph and
generalise this operation to hypergraphs. 
We establish lower bounds on the number of flat spectra of a Boolean function,
depending on internal structures,
with respect to the $\{I,H\}^n$ and $\{I,H,N\}^n$ sets of transforms.
We also construct a family of Boolean functions of degree higher than
two with a large number of flat spectra with respect to $\{I,H\}^n$,
and compute a lower bound on this number. 
The relationship between pivot orbits and equivalence classes of
error-correcting codes is then highlighted.
Finally, an enumeration of pivot orbits of various types of graphs is given, and
it is shown that the same technique can be used to classify codes.
\end{abstract}

\section{Introduction}

The {\em pivot} operation on a graph $G$ was used by Arratia, Bollob\'{a}s and
Sorkin~\cite{Arr:Int, Arr:Int1} 
to define the {\em interlace polynomial} $q(G,z)$,
as a variant of the Tutte and Tutte-Martin polynomials~\cite{Bou:Mart}. It was also described by 
Van den Nest~\cite{VdNe}, under the name of {\em edge-local complementation}.
In~\cite{RP:Int}, we related the interlace polynomials of a graph to the spectra of a 
quadratic Boolean function with respect to a strategic subset of local unitary transforms.
Our main motivation in doing this was to establish links between graph theory,
cryptography, coding theory, and quantum entanglement.  

Let the graph $G=(V,E)$, with vertex set, $V$, and edge set, $E$,
of order $n$ be represented by its $n \times n$ adjacency matrix,
$\Gamma$. Identify $G$
with a quadratic Boolean function $p(x_0,x_1,\ldots,x_{n-1})$, where
$p(\boldsymbol{x}) = \sum_{i < j} \Gamma_{ij}x_ix_j$~\cite{RP:BC1},
i.e., the term $x_ix_j$ occurs in $p(\boldsymbol{x})$ if and only if $ij \in E$.
This identification allows us to interpret $q(G,1)$ as the number of
{\em{flat spectra}} of $p(\boldsymbol{x})$
with respect to (w.r.t.) the set of transforms $\{I,H\}^n$.
In this paper we characterise the pivot operation using \emph{algebraic normal form} (ANF).
We also generalise pivot to
hypergraphs, and state the (necessary and sufficient) condition that a function
of degree higher than two must fulfil in order to allow such an operation.
Then we show how the pivot operation on a (hyper)graph can be written as a
transform from $\{I,H\}^n$ on the bipolar vector of the function associated to it.
We then prove that all (not necessarily all) flat spectra of a quadratic (general)
Boolean function, $p$,
w.r.t. $\{I,H\}^n$, can be realised
via a series of pivot operations on the graph (hypergraph) associated to $p$,
respectively.
We then construct a family of Boolean 
functions that have a large number of flat spectra w.r.t. $\{I,H\}^n$, and
compute this number.
We also study the pivot orbit of structures that include a
clique and develop lower bounds on the number of flat
spectra of a graph w.r.t. $\{I,H\}^n$ and $\{I,H,N\}^n$.
It is shown that orbits of bipartite graphs under the pivot operation correspond
to equivalence classes of binary linear codes, and that
all information sets of a code can be found by pivoting on its associated graph.
We also give an enumeration of pivot orbits of all graphs on up to 12 vertices, and
of all bipartite graphs on up to 13 vertices. 

To the best of our knowledge, the results mentioned above have not appeared in the literature before.

\section{Definitions and Notation}\label{DefNot}

Let $H = \frac{1}{\sqrt{2}} \begin{small}\left ( \begin{array}{rr}
1 & 1 \\
1 & -1
\end{array} \right )\end{small}$ be the Walsh-Hadamard kernel, 
$N = \frac{1}{\sqrt{2}} \begin{small}\left ( \begin{array}{rr}
1 & i \\
1 & -i
\end{array} \right )\end{small}$, where $i^2 = -1$,
be the Negahadamard kernel, and let $I$ the
$2\times 2$  identity matrix. Let $s = (s_0,s_1,....,s_{2^n-1})^t = (-1)^{p(\boldsymbol{x})}$, where 
$s_{\boldsymbol{i}} = (-1)^{p(\boldsymbol{i})}$ and $p(\boldsymbol{x}): \GF(2)^n \rightarrow \GF(2)$ is a Boolean function. With this notation, $p$ 
is {\em bent}~\cite{Rot:Bent} if $P=2^{-\frac{n}{2}}(\bigotimes_{i=0}^{n-1} H)
(-1)^{p(\boldsymbol{x})}$ has a {\em{flat}} spectrum, or, in other words, if
$P=(P_{\boldsymbol{k}}) \in {\mathbb{C}}^{2^n}$ is such that
$|P_{\boldsymbol{k}}| = 1$, $\forall \boldsymbol{k}\in \GF(2)^n$, where '$\otimes$' denotes the tensor
product of matrices, also known as the Kronecker product.
If the function is quadratic, we associate to it a simple undirected
graph, and in this case a flat spectrum is obtained if and only if $\Gamma$, the adjacency matrix 
of the graph, has maximum rank as a binary matrix.
In~\cite{RP:BC1}, we generalised this concept, considering not only the Walsh-Hadamard 
transform $\bigotimes_{i=0}^{n-1} H$, but the complete
set of $3^n$ unitary transforms $\{I,H,N\}^n$, comprising all transforms $U$ of the form
$U=\prod_{j \in {\boldsymbol{R_I}}} I_j \prod_{j \in {\boldsymbol{R_H}}} H_j \prod_{j \in {\boldsymbol{R_N}}} N_j,$
where the sets $\boldsymbol{R_I}, \boldsymbol{R_H}$ and $\boldsymbol{R_N}$ partition the set of vertices
$\{0,\ldots,n-1\}$, and  $H_j$, say, is short for
$I \otimes I \otimes \ldots \otimes I \otimes H \otimes I \otimes \ldots \otimes I$, 
with $H$ in the $j^{th}$ position.
For instance, if $n = 4$, ${\boldsymbol{R_I}} = \{1\}$, ${\boldsymbol{R_H}} = \{0,3\}$, and
${\boldsymbol{R_N}} = \{2\}$, then $U = H \otimes I \otimes N \otimes H$, where $U$ is a $16 \times 16$
unitary matrix. The \emph{orbit} of a Boolean function $p$ w.r.t. a set of transforms 
$T$ comprises all Boolean functions, $p'$, where $s'_{\bf{i}} =(-1)^{p'({\bf{i}})}$, and                where $s'$ can be obtained by the application of any $t\in T$ to $s = (-1)^{p(\boldsymbol{x})}$.

In~\cite{RP:BC1,RP:BC2}
we studied the number of flat spectra of a function w.r.t. $\{I,H,N\}^n$, 
or in other words the number of unitary transforms $U\in\{I,H,N\}^n$ such that 
${P_{U}}=(P_{U,{\boldsymbol{k}}})\in {\mathbb{C}}^{2^n}$ has $|P_{U,{\boldsymbol{k}}}|=1$, $\forall  \boldsymbol{k} 
\in \GF(2)^n$, where $(P_{U,{\boldsymbol{k}}})=U(-1)^{p({\boldsymbol{x}})}\enspace.\label{P}$
We also considered the number of flat spectra w.r.t. some subsets of $\{I,H,N\}^n$, namely $\{H,N\}^n$ (when ${\boldsymbol{R_I}}=\emptyset$) and $\{I,H\}^n$ (when ${\boldsymbol{R_N}}=\emptyset$).  
We also proved that a quadratic Boolean function will have a flat spectrum w.r.t. a transform $U\in\{I,H,N\}^n$ 
if and only if a certain modification of its adjacency matrix 
has maximal binary rank. 




As will be explained in the next section, the \emph{pivot orbit} of a (hyper)graph $G$ consists of all graphs obtained by the application of 
any sequence of \emph{pivot} operations to $G$. Similarly, the \emph{LC orbit} comprises
all graphs obtained by applying \emph{local complementations} to $G$.
In this paper, we will study the pivot orbits of (hyper)graphs, and the subsets of their 
LC orbits that are associated to the pivot transform.

There are two names for the pivot operation on graphs that are
currently in use in the
literature, namely {\em{pivot}} and {\em{edge-local complementation}} (ELC). The name
\textquotedblleft edge-local complementation" comes from Bouchet's original
definition of \textquotedblleft local
complementation on the edge" in \cite{Bou:Graph} and the name
\textquotedblleft edge-local complementation" has
been used recently by Van den Nest in \cite{VdNe}. The name \textquotedblleft pivot"
has a long history
with respect to Gaussian elimination and, in the context of graphs, would be the
operation of ELC on a bipartite graph. A few authors
\cite{Arr:Int,Arr:Int1,MonSar,riera:pivot} have, since Bouchet, extended the use of
\textquotedblleft pivot" to apply to all graphs, not just bipartite. We call the ELC operation,
\textquotedblleft pivot", in this paper, although we acknowledge
that \textquotedblleft edge-local complementation" is equally valid. Note, however, that in this paper
we further generalise to hypergraphs the applicability of pivot. 

\section{Pivot}


We recall the definition of two graph operations, \emph{local complementation} (LC),
also known as \emph{vertex neighbourhood complementation} (VNC), and \emph{pivot}, also known as {\em{edge-local complementation}} (ELC).

\begin{df}[\cite{Bou:Mart, Bou:Graph, Glynn:Graph, Hein:GrEnt, Fon:LC}]\label{LC}
Let $G=(V,E)$ be a graph and $i \in V$ be some vertex. $\mathcal{N}(i)$ denotes
the neighbourhood of $i$, i.e., the set of vertices adjacent to $i$.
The action of {\em{local complementation}} at vertex $i$, denoted LC($i$),
is defined as the graph transformation obtained by replacing $G[\mathcal{N}(i)]$, i.e., 
the subgraph induced on the neighbourhood of $i$, by its complement.
\end{df}

\begin{df}[\cite{Arr:Int,Arr:Int1,Bou:Graph}]\label{pivot}
Given a graph $G=(V,E)$ and an edge $ij \in E$, the action
of {\em{pivot}} on $ij$ is given by LC($i$)LC($j$)LC($i$), the action
of LC at vertex $i$, then vertex $j$, then vertex $i$ again.\footnote{Bouchet's original
definition of pivot, called \emph{complementation along an edge}~\cite{Bou:Graph},
includes a final swapping of the vertices $u$ and $v$.}
Note that the operations LC($j$)LC($i$)LC($j$) would give exactly the
same result.
\end{df}

\begin{figure}
 \centering
 \includegraphics[width=.40\linewidth]{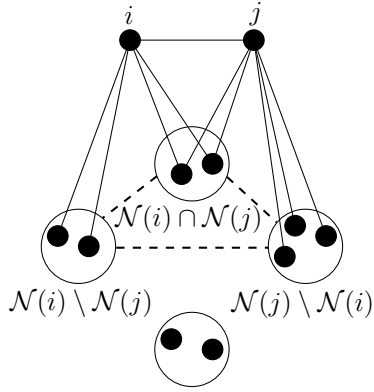}
 \caption{Visualisation of the Pivot Operation}\label{fig:pivot}
\end{figure}

\begin{df}[\cite{Arr:Int,Arr:Int1,Bou:Graph}]\label{pivot2}
Pivoting on $ij \in E$ of $G=(V,E)$ can equivalently be defined as follows.
Decompose $V \setminus \{i,j\}$ into four disjoint sets, as visualised in Fig.~\ref{fig:pivot},
\begin{itemize}
\item $\mathcal{N}(i) \setminus \mathcal{N}(j)$, the set of vertices adjacent to $i$, but not to $j$,
\item $\mathcal{N}(j) \setminus \mathcal{N}(i)$, the set of vertices adjacent to $j$, but not to $i$,
\item $\mathcal{N}(i) \cap \mathcal{N}(j)$, the set of vertices adjacent to both $i$ and $j$,
\item and the set of vertices adjacent to neither $i$ nor $j$.
\end{itemize}
For any pair of vertices $\{x,y\}$, where $x$ belongs to one of the first three classes listed above, 
and $y$ also belongs to one of the first three classes, but a different class than $x$, 
`toggle' the pair $\{x, y\}$, i.e., if $xy \in E$, delete the edge, 
and if $xy \not\in E$, add the edge $xy$ to $E$. Finally, swap the labels of vertices $i$ and $j$,
or, equivalently, swap the two sets $\mathcal{N}(i)$ and $\mathcal{N}(j)$.\footnote{In Bouchet's 
original definition of pivot, this swapping does not occur.}
\end{df}

Let the vertex $i$ of the graph $G=(V,E)$ correspond to the variable $x_i$ in $p(\boldsymbol{x})$, the
quadratic Boolean function associated to $G$.
As defined above, $\mathcal{N}(i)$ is the set of vertices that are adjacent to $i$.
We identify $\mathcal{N}(i)$ with the linear Boolean function
$\mathcal{N}_i = \sum_{k \in \mathcal{N}(i)} x_k$.
Thus $x_i \mathcal{N}_i$ is the quadratic Boolean function corresponding to all
edges incident on $i$. We can now redefine the pivot operation in terms of Boolean functions.

\begin{lem}\label{piv:ANF}
Let $p=x_ix_j+x_i{\cal{N}}_i+x_j{\cal{N}}_j+R$ be a quadratic Boolean function, where
${\cal{N}}_i$, ${\cal{N}}_j$, and $R$ are not functions of $x_i$ or $x_j$.
$p$ corresponds to the graph $G=(V,E)$, the linear function $\mathcal{N}_i$ corresponds to the 
neighbourhood of $i \in V$, $\mathcal{N}_j$ to the neighbourhood of $j \in V$, and the quadratic
function $R$ to all edges in $E$ that are incident on neither $i$ nor $j$.
The Boolean function corresponding to the graph obtained by pivoting on the edge $ij \in E$ is
\begin{equation}
\begin{split}
p_{iji} &= x_ix_j+x_i{\cal{N}}_j+x_j{\cal{N}}_i+{\cal{N}}_i{\cal{N}}_j+R\\
   &= p+(x_i+x_j)({\cal{N}}_i+{\cal{N}}_j)
+{\cal{N}}_i{\cal{N}}_j.
\end{split}
\end{equation}
Note that both $p$ and $p_{iji}$ can contain linear terms which may be ignored.
We consider $p$ and $p_{iji}$ to be equivalent, since the corresponding graphs
are equivalent up to pivot operations.
\end{lem}



\subsection{A Generalisation to Hypergraphs}

Let $p$ be a function of degree at least two. Let ${\cal{N}}_i$ now be the
Boolean function comprising all terms
which multiply $x_i$ in $p$, such that ${\cal{N}}_i$ is independent of $x_i$.
Note that there is no longer a 
restriction on the degree of ${\cal{N}}_i$.

\begin{df}\label{dependson}
For Boolean functions, $f$ and $g$,
we say that $g\in f$ or $g \notin f$,
when $f$ does or does not depend on $g$, respectively.
\end{df}

\begin{df}\label{term}
For Boolean functions, $f$ and $g$,
we say that $g \in_t f$ or $g \notin_t f$,
when $g$ is or is not a term in the algebraic
normal form of $f$, respectively.
\end{df}

\begin{df}\label{multiplying}
For Boolean functions, $f$ and $g$, we say that $g \in_m f$ and
$g \notin_m f$
when $g$ is or is not a {\em{multiplying term}} in $f$, respectively,
where $g$ is a multiplying term in $f$ iff $\exists r$ such that $gr \in_t f$.
\end{df}

\begin{df}\label{hyperpivot}
Let $p=x_ix_j+q(x_0,\ldots,x_{n-1})$ be a function of degree at least two
such that $x_ix_j \notin_m q$.
The function $p$ corresponds to the \emph{hypergraph} $G=(V,E)$, and $x_ix_j$ corresponds
to the edge $ij \in E$ of degree two. 
The Boolean function corresponding to the graph obtained by pivoting on $ij \in E$
is defined as
\begin{equation}
\begin{split}
p_{iji} &= x_ix_j+x_i{\cal{N}}_j+x_j{\cal{N}}_i+{\cal{N}}_i{\cal{N}}_j+R\\
 &= p+(x_i+x_j)({\cal{N}}_i+{\cal{N}}_j)+ {\cal{N}}_i{\cal{N}}_j,
\end{split}
\end{equation}
where $p=x_ix_j+x_i{\cal{N}}_i+x_j{\cal{N}}_j+R$ as before.
\end{df}

As a visualisation of pivot on hypergraphs, consider Fig.~\ref{fig:pivot}, where hyperedges
can be added anywhere, with the exception that no edge of degree higher than two can be incident on 
both $i$ and $j$. Due to (and equivalently to) the condition on $p$ in definition \ref{hyperpivot},
$\mathcal{N}_i$ and $\mathcal{N}_j$ are independent of both $x_i$ and $x_j$, and so the formula is 
well-defined. If we did not have this condition, the definition would have been ambiguous. When
$p$ is quadratic, and the vertices $i$ and $j$ of the corresponding graph are connected,
the condition is always fulfilled and the definition is consistent.
 
\begin{lem}\label{bip}
Let $G=(V,E)$ be a bipartite (hyper)graph. This means that
$V = X \cup Y$ such that none of the induced subgraphs $G[X]$ 
and $G[Y]$ contain any edges.
If we interpret $\boldsymbol{x}$ and $\boldsymbol{y}$ as vectors of variables,
representing the sets $X$ and $Y$, then $G$ corresponds to
a Boolean function $p=h(\boldsymbol{x}) \cdot g(\boldsymbol{y})$,
where $h(\boldsymbol{x})$ and $g(\boldsymbol{y})$ are vectors of Boolean functions of any degree. 
After pivoting on any permissible edge of $G$, the resulting (hyper)graph always remains bipartite.
Moreover, the sizes of the two partitions will not change under pivot operations.
\end{lem}


\begin{proof}
It follows from the definition of a bipartite (hyper)graph that
for any edge $ij \in E$, $i$ and $j$ have no common neighbours, and the subgraphs of $G$ induced 
on $\mathcal{N}(i)$ and $\mathcal{N}(j)$ contain no edges. 
It follows from lemma~\ref{piv:ANF} 
for graphs and definition~\ref{hyperpivot} for hypergraphs that the (hyper)graph obtained
by pivoting is also bipartite.
\end{proof}


\subsection{Pivot in Spectral Terms}

In~\cite{RP:BC1}, we proved that local complementations on a graph can be realised via the application 
of successive negahadamard ($N$) transforms on the bipolar vector, $s=(-1)^p$, of the associated function $p$.
We here show that pivot operations on a (hyper)graph also correspond to certain transformations from
the set $\{I,H,N\}^n$.

Let $m : \GF(2)^n \rightarrow \GF(2)$.
In the following, we shall embed the output of $m$ in the complex 
numbers by 
the operation $[m] \in {\mathbb{C}}$, where $[0] = 0$, and $[1] = 1$.\footnote{Note that $[1+1]=[0]=0$, while
$[1]+[1]=1+1=2$}

Let   $s =[m({\bf{x}})](-1)^{p({\bf{x}})}$ be dependent on all binary variables
$x_i,\ 0\leq i\leq n-1$, where 
$m = \prod_{k=0}^{u-1} h_k$ and the $h_k$ are Boolean functions in $n$
variables,\begin{footnote}{Such a factorisation of $m$ is not necessarily unique.}\end{footnote}
and $p$ is a Boolean
function of degree less or equal than two. In the sequel, expressions of the form $s = c[m](-1)^p$, with
$c \in {\mathbb{C}}$, shall always be written as $s = [m](-1)^p$, i.e. we ignore
normalisation coefficients. For an index $i$, we write $m=rv$, where all the 
terms in $v=\prod_{k\in V} h_k$, for
some $V\subseteq\{0,\ldots,u-1\}$,  depend on $x_i$, and $r$ does
not depend on $x_i$. We denote $p_a=p|_{x_i=a},\ m_a=m|_{x_i=a},\ v_a=v|_{x_i=a},$ for
$a \in \GF(2)$.
From the conditions above, and by results of \cite{PR:QE}, we get the
following theorems.

\begin{theorem}\label{general}
Let $s=[m](-1)^p$. Then 
\begin{equation}
H_is= [r(v_0+v_1)](-1)^{p_0+v_1(p_0+p_1+x_i)} + 2[rv_0v_1(p_0+p_1+x_i+1)](-1)^{p_0}.
\end{equation}
\end{theorem}
\begin{proof}
$ s = [m](-1)^{p} = [(1 + x_i)m_0](-1)^{p_0} + [x_im_1](-1)^{p_1} $.
Applying $H_i$ gives,
\begin{equation}\label{eqa1}
\begin{split}
s' &= [1 + x_i]([m_0](-1)^{p_0} + [m_1](-1)^{p_1}) + [x_i]([m_0](-1)^{p_0} - [m_1](-1)^{p_1}) \\
  &= [1 + x_i]([m_0(p_0 + 1)] + [m_1(p_1 + 1)] - [m_0p_0] - [m_1p_1])\\
  &\quad + [x_i]([m_0(p_0 + 1)] - [m_1(p_1 + 1)] - [m_0p_0] + [m_1p_1])
\end{split}
\end{equation}
By applying the following identity to (\ref{eqa1}),
for Boolean functions $A_0,A_1,B_0,B_1$,
\[
\begin{split}
[A_0] + [A_1] + [B_0] + [B_1] & = [A_0 + A_1 + B_0 + B_1](-1)^{A_0A_1 + B_0B_1 + B_0 + B_1} \\
 &\quad + 2[(A_0 + A_1 + B_0 + B_1)(A_0A_1 + B_0B_1)](-1)^{A_0 + 1},
\end{split}
\]
we obtain, after a bit more manipulation, the theorem.\begin{footnote}{
Theorem \ref{general} and its proof relate to theorem 17 of \cite{PR:QE}.
However, we have included a new proof as the proof
of theorem 17 was incorrect in \cite{PR:QE}. We have also simplified the statement
of the theorem.
}\end{footnote}
\end{proof}

\begin{theorem}[theorem 18 of \cite{PR:QE}]\label{notinm}
Let $s=[m](-1)^p$. If $x_i\notin m$, then 
\begin{equation}
H_is=[m\cdot(p_0+p_1+x_i+1)](-1)^{p_0}.
\end{equation}
\end{theorem}

\begin{theorem}[theorem 20 of \cite{PR:QE}]\label{inm}
Let $s=[m](-1)^p$. If $x_i\in m$ and if there exists a factorisation of $v$
such that all $h_k \in_m v$ are linearly dependent on $x_i$, then 
\begin{equation}
H_is=[r\cdot(v_0+v_1)](-1)^{p_0+h_{z,1}(p_0+p_1+x_i)},
\end{equation}
where $h_{z,1}=h_{z}|_{x_i=1}$ and $v_0+v_1=\prod_{k\neq j} (h_j+h_k+1)$, with $h_z$ and $h_j$ chosen arbitrarily 
among the divisors of $v$. 
\end{theorem}

\begin{remark}
Typically we will choose $z=j$.
\end{remark}

\begin{theorem}\label{pivspec}
Let $p$ be a Boolean function that fulfils the condition of definition~\ref{hyperpivot}. Then any 
(hyper)graph obtained by pivoting on the (hyper)graph associated to $p$ corresponds to some member 
of the set of $\{I,H\}^n$ transforms of $p$. Concretely, if $p_{iji}$ is the function obtained by pivoting
on the edge $ij$ of the (hyper)graph associated with $p$, then  
$(-1)^{p_{iji}} = (H_i\cdot H_j)(-1)^p$.
\end{theorem}

\begin{proof} 
Let $p=x_ix_j+x_i{\cal{N}}_i+x_j{\cal{N}}_j+R$, and let $s=(-1)^p$. Then, by theorem~\ref{notinm},
\begin{equation}
s'=H_is=[x_j+{\cal{N}}_i+x_i+1](-1)^{x_j{\cal{N}}_j+R}.
\end{equation}
Now, applying theorem~\ref{inm}, we get 
\begin{equation}
s''=H_js'=
1\cdot (-1)^{R+({\cal{N}}_i+x_i)({\cal{N}}_j+x_j)}=(-1)^{x_ix_j+x_i{\cal{N}}_j
+x_j{\cal{N}}_i+{\cal{N}}_i{\cal{N}}_j+R},
\end{equation}
which is what we wanted. By the condition on $p$, ${\cal{N}}_i$ does not depend 
on $x_j$, which ensures that the conditions on $m$ necessary to apply
theorem~\ref{inm} are fulfilled.
\end{proof}

\begin{cor}\label{corhyperpivot}
Let $p$ be a Boolean function of any degree that satisfies the conditions of 
definition~\ref{hyperpivot}. Then $p$ has a flat spectrum with respect to the
transform $U = H_i \cdot H_j$. 
\end{cor}

\begin{theorem}\label{OnlyPivot}
Each of the flat spectra of $p$ with respect to the set of transforms
$\{H_i \cdot H_j \mid| i,j \in \mathbb{Z}_n, i \ne j\}$, identifies
an edge $ij$ in the hypergraph associated with $p$, and therefore can be obtained
by pivoting on the hypergraph associated with $p$ at the edge $ij$.
\end{theorem}
\begin{proof}
We will show that, given some arbitrary spectrum, $H_i(-1)^p$, the only way one
can obtain a flat spectrum, $H_j \cdot H_i(-1)^p$, $i \ne j$,
is when
$$ x_ix_j \in_t p, \hspace{5mm} x_ix_j \notin_m p - x_ix_j. $$
In such a case, theorem \ref{pivspec} states that
$H_j \cdot H_i(-1)^p$ always corresponds to a pivot operation on the
hypergraph associated to $p$ at the edge $ij$.

From theorem \ref{notinm}, for arbitrary $i$,
$$ H_i(-1)^p = [p_0 + p_1 + x_i + 1](-1)^{p_0} = [{\tilde{m}}](-1)^{p_0}, $$
for some ${\tilde{m}}$.
In order that $H_j \cdot H_i(-1)^p = (-1)^{p'}$, for some $p'$,
we must transform the factor,
$[p_0 + p_1 + x_i + 1]$, back to $1$. This is trivially possible if $j = i$, but
the theorem excludes the case where $i = j$.
Let
$$ H_j \cdot H_i(-1)^p = [m'](-1)^{p'}, $$
for some $m'$ and $p'$,
where $i$ and $j$ are arbitrary, $i \ne j$. We wish to choose $j$ such that $m' = 1$.
There are three possible scenarios:
\begin{itemize}
	\item $x_j \notin {\tilde{m}}$: In this case, from theorem \ref{notinm},
	$(p_0 + p_1 + x_i + 1) \in_m m'$ so $m' \neq 1$.
	\item $x_j \in {\tilde{m}}$: There are three subcases. Let
	$$ p_{00} = p_0|_{x_j=0}, \hspace{3mm} p_{10} = p_1|_{x_j=0},
	\hspace{3mm} p_{01} = p_0|_{x_j=1}, \hspace{3mm} p_{11} = p_1|_{x_j=1}. $$
	Considering theorem \ref{general} acting on $[{\tilde{m}}](-1)^{p_0}$,
	then $m'$ can be $1$ iff one or
	more of the following three conditions are met:
\begin{eqnarray}
v_0 + v_1 = 1, \quad v_0v_1(p_{00}+p_{01}+x_j+1) &= 0\label{cond1}\\
v_0 + v_1 = 0, \quad v_0v_1(p_{00}+p_{01}+x_j+1) &= 1\label{cond2}\\
v_0 + v_1 = 1, \quad v_0v_1(p_{00}+p_{01}+x_j+1) &= 1\label{cond3}
\end{eqnarray}
As, in this case, $v = p_0 + p_1 + x_i + 1$,
$v_0 = v|_{x_j=0}$, $v_1 = v|_{x_j=1}$, then we observe that
$v_0 + v_1 = p_{00} + p_{10} + p_{01} + p_{11}$ and 
$v_0v_1 = (p_{00} + p_{10} + x_i + 1)(p_{01} + p_{11} + x_i + 1)$.
The three subcases for $x_j \in {\tilde{m}}$ are:
	\begin{itemize}
		\item $x_ix_j \notin_m p$:
		   In this case $v_0 + v_1 = 0$ so (\ref{cond1}) and (\ref{cond3}) are
		   not satisfied. Furthermore, $\deg(v_0v_1(p_{00}+p_{01}+x_j+1)) > 0$ as 
		   $v_0v_1 =p_{00}+p_{10}+x_i+1=p_{01} +
 p_{11} + x_i + 1$, so $x_i\in v_0v_1, x_j\notin v_0v_1$, and
 $x_i\notin p_{00}+p_{01}+x_j+1$, $x_j\in p_{00}+p_{01}+x_j+1$, 
so
		   (\ref{cond2}) is not satisfied.
		\item $x_ix_j \in_t p$, $x_ix_j \notin_m p - x_ix_j$:
		   In this case $v_0 + v_1 = 1$. Moreover,
		   $p_{00} + p_{10} = p_{01} + p_{11} + 1$, so $v_0v_1 = 0$.
		   Therefore (\ref{cond1}) is satisfied.	
		\item $x_ix_j \in_t p$, $x_ix_j \in_m p - x_ix_j$:
		   In this case $\deg(v_0+v_1) > 0$ so none of
		   (\ref{cond1}),(\ref{cond2}), or (\ref{cond3}) are satisfied.
	\end{itemize}
\end{itemize}
From the above analysis, $m' = 1$ iff $x_ix_j \in_t p$ and
$x_ix_j \notin_m p - x_ix_j$. This is precisely the condition required to ensure
pivot at the edge $ij$ on the hypergraph associated to $p$, as stated by
definition \ref{hyperpivot}.
\end{proof}

\begin{theorem}
Let $p$ be a quadratic Boolean function over $n$ variables. Then all
flat spectra of $p$ with respect to transforms fron the set $\{I,H\}^n$,
other than the identity, can be obtained via a sequence of pivot operations on the
graph associated to $p$.
\label{quadpiv}
\end{theorem}
\begin{proof}
Consider the following hypothesis:

Let ${\bf{X}}$ be a fixed subset  of
$\{0,1,\ldots,n-1\}$, where $|{\bf{X}}| > 2$. Let
$U = \prod_{i \in {\bf{X}}} H_i$. Then it is possible for
$U(-1)^p$ to be flat, and for
$U'(-1)^p$ not to be flat $\forall\ U'$ satisfying
$U' = \prod_{i \in {\bf{Z}}} H_i$, where ${\bf{Z}} \subset {\bf{X}}$ and
${\bf{Z}} \neq \emptyset$.

The theorem is proved if the hypothesis is proven false, as $H(-1)^p$ is never flat. 
If the hypothesis is true for some ${\bf{X}}$, then
$\nexists\ i,j  \in {\bf{X}}$ such that $H_i \cdot H_j(-1)^p$ is flat. We know, from
theorem \ref{pivspec} that, therefore, the set of vertices, ${\bf{X}}$,
forms an independent set\begin{footnote}{If ${\bf{X}}$ is not an independent set, then there is an edge
 $ij$ between vertices of ${\bf{X}}$, and thereby we can pivot on it, and $H_i \cdot H_j(-1)^p$ is flat.}\end{footnote} in the graph, $G$, associated to $p$. But
$U(-1)^p$ cannot be flat if ${\bf{X}}$ is an independent set in $G$ as,
applying $H$ to $(-1)^p$ at all index positions in ${\bf{X}}$ requires
$|{\bf{X}}|$ invocations of theorem \ref{notinm}, each of which contributes
a new linear factor to $m$. Therefore the final $m$ cannot be $1$ and the
hypothesis is false. But, for $|{\bf{X}}| = 2$, we know from theorem
\ref{OnlyPivot} that all flat spectra are obtained via pivot operations. It is
trivial to show that $U(-1)^p$ is never flat if $|{\bf{X}}| = 1$.
\end{proof}

\begin{lem}
Let $p$ be a Boolean function of any degree over $n$ variables. Then there may exist
flat spectra of $p$ with respect to transforms fron the set $\{I,H\}^n$,
other than the identity, that cannot be obtained via a sequence of pivot
operations on the hypergraph associated to $p$.
\label{genpiv}
\end{lem}
\begin{proof}
By example, the Boolean function,
\[
\begin{split}
p(\boldsymbol{x}) &= x_0x_1x_2+x_0x_1x_3+x_0x_1x_5+x_0x_2x_4+x_0x_2x_5+x_0x_3x_4+x_0x_3x_5 \\
 &\quad +x_0x_4x_5+x_1x_2x_3+x_1x_2x_4+x_1x_2x_5+x_1x_3x_4+x_1x_4x_5+x_2x_3x_4\\
 &\quad +x_2x_3x_5+x_3x_4x_5,
\end{split}
\]
has two flat spectra w.r.t. the set $\{I,H\}^6$. Apart from the identity
transform, $(-1)^p$ is also flat w.r.t.
$H \otimes H \otimes H \otimes H \otimes H \otimes H$. Such a flat spectrum
cannot be obtained via a series of pivot operations as $p$ does not contain
any quadratic terms.
\end{proof}
\begin{remark}
The example used in the proof of lemma \ref{genpiv} was taken from
an interesting catalogue of {\em{homogeneous bent functions}}, as provided by
\cite{Char:HomBent}.
\end{remark}

Let $p : \GF(2)^n \rightarrow {\mathbb Z}_4$,
$m : \GF(2)^n \rightarrow \GF(2)$, and let $s \in {\mathbb C}^{2^n}$ be such 
that $s = (s_0,s_1,\ldots,s_{2^n-1})^t = [m(\boldsymbol{x})]i^{p(\boldsymbol{x})}$, where 
$s_{\boldsymbol{j}} = [m(\boldsymbol{j})]i^{p(\boldsymbol{j})}$. Sometimes, for brevity, we write the above as 
$s = [m]i^p$, when it is clear from the context what we mean.  Let $m_a : \GF(2)^n \rightarrow \GF(2)$ 
represent $m_a = m|_{x_j = a}$. Similarly, let $p_a : \GF(2)^n \rightarrow {\mathbb Z}_4$ represent 
$p_a = p|_{x_j = a}$.

\begin{theorem}\label{NAPFext}
 Let $m:\GF(2)^n\rightarrow \GF(2)$ and $p:\GF(2)^n\rightarrow {\mathbb Z}_4$.
 Then,
\begin{equation}
N_j[m]i^p=\frac{1}{\sqrt{2}}([m_0]i^{p_0}+[m_1]i^{p_1+2x_j+1}).
\end{equation}
\end{theorem}

\begin{proof} 
Without loss of generality, we set $j=n-1$. Then, we can write the complex vector
$[m]i^p$ (seen as a $2^n \times 1$ matrix) as 
\[
[m]i^p = \left(\begin{array}{c}
[m_0] i^{p_0}\\
{[m_1]} i^{p_1}
\end{array}\right),
\]
where $[m_0] i^{p_0}$ and $[m_1] i^{p_1}$ are $2^{n-1}\times 1$ complex matrices.
Then,
\[
\begin{split}
 N_{n-1}[m]i^p  &= 
 \left(\begin{array}{cccc}
 N & \boldsymbol{0} & \ldots & \boldsymbol{0}\\
 \boldsymbol{0} & N & \ldots & \boldsymbol{0}\\
 \vdots & \vdots & \ddots & \vdots\\
 \boldsymbol{0} & \boldsymbol{0} & \ldots & N 
 \end{array}\right)
 \left(\begin{array}{c}
 [m_0]i^{p_0}\\
 {[m_1]}i^{p_1} 
 \end{array}\right)\\
&=
 \left(\begin{array}{c}
 [m(0,\ldots,0)]i^{p(0,\ldots,0)}+i[m(0,\ldots,1)]i^{p(0,\ldots,1)}\\
 {[m(0,\ldots,0)]}i^{p(0,\ldots,0)}-i[m(0,\ldots,1)]i^{p(0,\ldots,1)}\\
 {[m(0,\ldots,1,0)]}i^{p(0,\ldots,1,0)}+i[m(0,\ldots,1,1)]i^{p(0,\ldots,1,1)}\\
 {[m(0,\ldots,1,0)]}i^{p(0,\ldots,1,0)}-  i[m(0,\ldots,1,1)]i^{p(0,\ldots,1,1)}\\
 \vdots\\
 {[m(1,\ldots,1,0)]}i^{p(1,\ldots,1,0)}+i[m(1,\ldots,1)]i^{p(1,\ldots,1,1)}\\
 {[m(1,\ldots,1,0)]}i^{p(1,\ldots,1,0)}-i[m(1,\ldots,1,1)]i^{p(1,\ldots,1,1)} \end{array} \right)
 \\
  &= [m_0]i^{p_0}+[m_1]i^{p_1+2x_{n-1}+1}.
\end{split}
\]
\end{proof}

In~\cite{RP:BC1}, we proved that local complementation can be realised via the application of 
successive $N$s on the bipolar vector of the function,
$s=(-1)^p$. Let $D$ be the set of (unitary) 
diagonal or anti-diagonal $2\times 2$ matrices.
Define $\delta,\,\gamma\in \{D\}^n$ as 
$\delta=\frac{\sqrt{2}}{1+i}\displaystyle\prod_{k=l,j}
\begin{small}\left(\begin{array}{cc}
1 & 0\\
0 & i
\end{array}\right)_{\!\!k}
\end{small}$ and $\gamma=-\begin{small}\displaystyle\prod_{k=l,j}
\left(\begin{array}{cc}
 0&-1\\
 1&0\end{array}\right)_{\!\!k}
\end{small}$.  Then, 

\begin{theorem} 
Let $p$ be a function that fulfils the condition of definition~\ref{hyperpivot}. Then the 
local complementation of its associated (hyper)graph, seen as a weighted (hyper)graph, lies in the orbit of $p$
w.r.t. $\{I,H,N\}^n$ to within a 
post-multiplication by a tensor product of members of $D$. Concretely, if $p_l,\, p_{jl}$,
and $p_{ljl}$ are
the functions obtained by local complementations on the vertices $l$, $j$,
then $l$ again, of the (hyper)graph associated with $p$, then 
\begin{equation}
\begin{array}{l}   i^{p_l} = \delta N_l(-1)^p, \\
        i^{p_{jl}} = \delta N_j\delta N_l(-1)^p,\\
        (-1)^{p_{ljl}} = \gamma\delta N_l\delta N_j\delta N_l(-1)^p.
\end{array}
\end{equation}
\end{theorem}

\begin{remark}
We do not distinguish between $l$ and $j$, so one can obtain the hypergraphs
associated to 
the functions, $p_l$, $p_j$, $p_{jl}$, $p_{lj}$, and $p_{ljl}$, via local complementation.
Note that $p_{jlj} = p_{ljl}$.
\end{remark}


\begin{proof}  
Let $p=x_lx_j+x_l{\cal{N}}_l+x_j{\cal{N}}_j+R$, and $s=(-1)^p$.
Let ${\cal{N}}_l=\sum_{r=0}^\rho u_r$, and 
${\cal{N}}_j=\sum_{t=0}^\tau v_t$, (note that they are not necessarily linear). Then, 
applying theorem~\ref{NAPFext}\footnote{One can lift the Boolean function $p$ to a function
$q:\GF(2)^n \rightarrow {\mathbb Z}_4$, with $q(\boldsymbol{x})=2p(\boldsymbol{x})$.}
(or by the results on \cite{RP:BC1}), $N_ls=\frac{1+i}{\sqrt{2}}i^{p'}$, where 
$p':\GF(2)^n\rightarrow {\mathbb Z}_4$, 
with explicit formula\footnote{We denote as $\lambda_0\phi_0 + \lambda_1\phi_1$ or,
more generally, as 
$\sum\lambda_i\phi_i$, with $\lambda_i\in {\mathbb Z}_4$ and $\phi_i$ Boolean
functions, the result of embedding the output of the $\phi_i$'s into ${\mathbb Z}_4$, multiplying them by a scalar 
$\lambda_i\in {\mathbb Z}_4$, and then adding the result mod 4. For instance, for $x_0 = x_1 = 1$, we have $2[x_0+x_1] +3x_1 + 2 = 1$.}
\begin{equation}
p'=2\left(p(x)+x_j\sum_{r=0}^\rho u_r+ 
\sum_{r\neq s}u_ru_s\right)
+3\left(x_l+x_j+\sum_{r=0}^\rho u_r\right).
\end{equation}
Applying $\delta$ to $N_ls$, we get $s'=\delta
N_ls=i^{p_l}$, where 
\begin{equation}
p_l=2\left(p(x)+x_j\sum_{r=0}^\rho u_r+
\sum_{r\neq s}u_ru_s\right)
+3\sum_{r=0}^\rho u_r.
\end{equation} 
This is the result of the action of LC$(l)$.
Now we apply LC$(j)$; that is, we first apply $N_j$ to $s'$. By theorem~\ref{NAPFext}, the result is 
$N_js'=\frac{1+i}{\sqrt{2}}i^{p''}$, where $p'':\GF(2)^n\rightarrow {\mathbb Z}_4$, with explicit formula 
\begin{equation}
\begin{split}
p'' &= 2\left(x_lx_j+x_l\sum_{t=0}^\tau v_t+x_j\left(\sum_{r=0}^\rho u_r+
\sum_{t=0}^\tau v_t\right)\right.\\
&+ \left. \sum_{t\neq u}v_tv_u+\sum_{r,t}u_rv_t+\sum_{r=0}^\rho u_r+R\right)\\
&+ 3(x_l+x_j+\sum_{t=0}^\tau v_t)
\end{split}
\end{equation}
Then we apply $\delta$ to $N_js'$ to get $s''=\delta
N_js'=i^{p_{lj}}$, where 
\begin{equation}
\begin{split}
p_{lj} &= 2\left(x_lx_j+x_l\sum_{t=0}^\tau v_t+x_j\left(\sum_{r=0}^\rho u_r+
\sum_{t=0}^\tau v_t\right)\right.\\
&+ \left. \sum_{t\neq u}v_tv_u+\sum_{r,t}u_rv_t+\sum_{r=0}^\rho u_r+R\right) + 3\sum_{t=0}^\tau v_t
\end{split}
\end{equation}
Now we apply LC$(l)$ again; that is, we first apply $N_l$ to $s''$. By theorem~\ref{NAPFext},
the result is $N_ls''=\frac{1+i}{\sqrt{2}}i^{p'''}$, 
where $p''':\GF(2)^n\rightarrow {\mathbb Z}_4$, with explicit formula 
\begin{equation}
\begin{split}
p''' &= 2\left(x_lx_j+x_l\sum_{t=0}^\tau v_t+x_j\sum_{r=0}^\rho u_r\right.\\
&+ \left. \sum_{r,t}u_rv_t+\sum_{r=0}^\rho u_r+\sum_{t=0}^\tau v_t+R\right)
+\displaystyle 3(x_l+x_j) 
\end{split}
\end{equation}
Then we apply $\delta$ to $N_ls''$ to get $s'''=\delta N_ls''=(-1)^{p_{ljl}'}$, where 
\begin{equation}
p_{ljl}'= x_lx_j+x_l\sum_{t=0}^\tau v_t+x_j\sum_{r=0}^\rho u_r
+\sum_{r,t}u_rv_t+\sum_{r=0}^\rho u_r+\sum_{t=0}^\tau v_t+R
\end{equation}
If we now apply $\gamma$ to $s'''$, we get 
\begin{equation}
p_{ljl}=\displaystyle x_lx_j+x_l\sum_{t=0}^\tau v_t+x_j\sum_{r=0}^\rho u_r
+\sum_{r,t}u_rv_t+R,
\end{equation}
which is, by definition \ref{hyperpivot}, the formula for pivot on the hypergraph
associated to $p$.   
Note that this  gives as well an alternative proof of theorem \ref{pivspec}: Let 
$d=\begin{small}\left(\begin{array}{cc}
1&0\\
0&i\end{array}\right)\end{small}$, 
and let $d'=\begin{small}\left(\begin{array}{cc}
0&-1\\
1&0\end{array}\right)\end{small}$. 
We see that we have applied:
\begin{itemize}
\item In position $l$: $d'dNddN=d'\begin{small}\left(\begin{array}{cc}
 0&1\\
-1&0\end{array}\right)\end{small}H=H$
\item In position $j$: $\frac{-1}{e^{3\pi i/4}}d'ddNd=(-1)d'\begin{small}\left(\begin{array}{cc}
0&-1\\
1&0\end{array}\right)\end{small}H=H$
\item Remaining positions: $I$
\end{itemize}
\end{proof}


\section{Number of Flat Spectra w.r.t. $\{I,H\}^n$}

In this section, we first study the behaviour of a graph that contains a certain inner structure, namely a {\em
clique}. Then we give bounds on the number of flat spectra of
graphs based on their subgraphs, specifically in the case where some of the subgraphs are cliques.

In order to provide some context for the results, we first state the results of some computer experiments. Table \ref{tab:flatIH} shows the expected number of flat spectra w.r.t. $\{I,H\}^n$ for a random Boolean function, and for a random Boolean function of degree $\le 2$. Table \ref{tab:flatIH} demonstrates,
empirically\begin{footnote}{Exhaustive search for random, $n \le 4$, and random quad., 
$n \le 7$, otherwise 100000 samples were taken.}\end{footnote} that, for $n$ large enough, the
expected number of flat spectra w.r.t. $\{I,H\}^n$ for a random Boolean function, and
for a random Boolean function of degree $\le 2$, respectively, is $1.0$ and approximately
$2^{n-2}$, respectively. The structures and  constructions considered in this section
will be seen to produce (hyper)graphs with relatively high numbers of flat spectra w.r.t.
$\{I,H\}^n$, in comparison to the average.
	     
\begin{table}
\centering
\caption{Average number of flat spectra w.r.t. $\{I,H\}^n$}
\label{tab:flatIH}
\begin{tabular}{cccccccccc}
\toprule
& \multicolumn{8}{c}{$n$}\\
\cmidrule(l){2-9}
 & 2 & 3 & 4 & 5 & 6 & 7 & 8 & 9\\
\midrule
 random & 1.500 & 1.750 & 1.390 & 1.039 & 1.000 & 1.000 & 1.000 & 1.000 \\
random quad  & 1.500 & 2.500 & 4.438 & 8.188 & 15.486 & 29.726 & 57.918 & 113.227 \\
\bottomrule
\end{tabular}
\end{table}

\subsection{Cliques}~\label{sec:clique}

The \emph{complete graph} or {\em clique} on $n$ vertices corresponds to the Boolean function
$p = \sum_{0\leq i<j\leq n-1} x_i x_j$.

\begin{lem}[\cite{RP:BC2}]\label{clique}
The Boolean function corresponding to the complete graph on $n$ vertices 
has $2^{n-1}$ flat spectra w.r.t. $\{I,H\}^n$, and
maximises over the set of Boolean functions of $n$ variables
the number of flat spectra w.r.t. $\{I,H\}^n$. 
\end{lem}

We now study the behavior of a graph that contains a \emph{clique}, i.e., a complete subgraph.
We consider three cases, depending on the positions of the vertices $a$ and $b$, where we pivot on the edge 
$ab$. Let $C_r$ be the clique on $r$ vertices
contained in the graph. We denote by $\mathcal{N}(a)$ and $\mathcal{N}(b)$ the neighbourhoods of $a$ and $b$.

\begin{itemize}
\item $a, b\in C_r$: The clique remains invariant.
\item $a \in C_r$, $b\notin C_r$: Let $m$ be the number of variables of $C_r$ that are in 
$\mathcal{N}(a) \cap \mathcal{N}(b)$.
Then, $C_r$ splits and we get
the cliques $C_{r-m},\ C_{m+2}$, connected just by $b$. 
Moreover $a \notin C_{r-m}$, $b \in C_{r-m}$, and $a, b\in C_{m+2}$.

\begin{itemize}
\item Particular case: Two connected cliques: $a \in C_{r_a}$, $b \notin C_{r_a}$, and $b\in C_{r_b}$.
  Let $m_a$ be the number of vertices of the clique $C_{r_a}$ that are in $\mathcal{N}(a) \cap \mathcal{N}(b)$,
  and $m_b$ the number of vertices of the clique $C_{r_b}$ that are in $\mathcal{N}(a) \cap \mathcal{N}(b)$. Then,
  both cliques split and we get the cliques $C_{r_a-m_a}$, $C_{m_a+2}$, $C_{r_b-m_b}$, and $C_{m_b+2}$. 
\end{itemize}

\item $a, b \notin C_r$: In this case, $C_r$ remains invariant, independently of
whether $a$ or $b$ are connected to it or not.
\end{itemize}

\subsection{Bounds on the Number of Flat Spectra}

We give lower bounds on the number of flat spectra w.r.t. $\{I,H\}^n$ and $\{I,H,N\}^n$ 
depending on internal structures:

\begin{lem} 
Consider an unconnected graph $G$, composed of two connected components, $G_1$ and $G_2$.
The number of flat spectra of $G$ w.r.t. $\{I,H\}^n$, 
$K_{IH}$, has as lower bound: $K_{IH}(G) \geq K_{IH}(G_1)\cdot K_{IH}(G_2)$ 
\end{lem}
 
\begin{cor} If we decompose an unconnected graph into connected its components $G_1, \ldots, G_t$, then 
$K_{IH}(G)\geq \prod_{i=1}^t K_{IH}(G_i)$. For instance, if we can decompose the graph into cliques 
$C_{r_1}, \ldots, C_{r_t}$, of respective sizes $r_1,\ldots, r_t$, then 
$K_{IH}(G)\geq \prod_{i=1}^t 2^{n_i-1}$.
\end{cor}

\begin{lem} 
Consider the number of flat spectra w.r.t. $\{I,H,N\}^n$. If we decompose an unconnected graph into 
connected components $G_1, \ldots, G_t$, then we have that $K_{IHN}(G)\geq \prod_{i=1}^t K_{IHN}(G_i)$.
\end{lem}

\begin{cor} 
The maximum clique size, $n_x$, of any member of the pivot orbit of $G$ is upper-bounded by 
$n_x \le \lfloor \log_2(K_{IH}) \rfloor$. 
\end{cor}

\section{A Construction of Boolean Functions with High Number of Flat Spectra}

We now design a family of Boolean functions in $n$ variables of 
degree less than or equal to $\max\{t,2\}$, where $0 \leq t \leq n-1$, 
whose members have a large number of flat spectra w.r.t. $\{I,H\}^n$.
Let 
\begin{equation}
f^{n,t}=\sum_{i=0}^{t-1}\sum_{j=t}^{n-1}x_ix_j + 
\sum_{i=t}^{n-2}\sum_{j=i+1}^{n-1}x_ix_j + a(x_0,x_1,\ldots,x_{n-1}),
\end{equation}
where $\deg(a) \le 1$. We then define the family ${\cal F}^{n,t}$,
\begin{equation}
{\cal F}^{n,t}=\{f^{n,t}+ h(x_0,x_1,\ldots,x_{t-1}) \},
\end{equation}
where $h$ is an arbitrary Boolean function of $t$ variables.

\begin{conjecture}
Let $f\in {\cal F}^{n,t}$. Then the pivot orbit of $f$ is contained in
$\displaystyle\bigcup_{k=0}^{n-1} {\cal F}^{n,k}$.
\end{conjecture}





\begin{theorem}\label{family}
Let $f\in {\cal F}^{n,t}$. Then the number of flat spectra of $f$ w.r.t. $\{I,H\}^n$ is at least 
$(t+1)2^{n-t-1}$, where the bound is tight if $f$ has degree $t$.
\end{theorem}

\begin{proof} 
Let $f\in {\cal F}^{n,t}$. Then it fulfils the condition of definition~\ref{hyperpivot} for every
edge $ij$ such that $t \leq i,j \leq n$. We showed in Section~\ref{sec:clique} that pivoting on 
any of these edges leaves the clique invariant. This means that the number of flat spectra of $f$ 
will be at least the number of times we can pivot on the clique on the last $n-t$ variables times 
the number of times we can pivot on the complete bipartite graph $\sum_{i=0}^{t-1}\sum_{j=t}^{n-1}x_ix_j$ 
(not counting repetitions), plus the identity transform. The number of times we can pivot on
the clique of the hypergraph is the same as the number of times we can pivot on a clique of size $n-t$. 
By lemma~\ref{clique}, this number is $2^{n-t-1}$. We can pivot on each edge of the complete bipartite 
graph, but note that the pivoting now changes the graph, so a new pivot may not be possible 
(depending on $h(x_0,\ldots, x_{t-1})$). Avoiding repetitions, that makes one pivot for every vertex on 
the first $t$ variables, plus the identity transform. In total we get the lower bound $(t+1)2^{n-t-1}$.

Let $f\in {\cal F}^{n,t}$ such that its degree is $t$. Take $h(x_0,x_1,\ldots,x_{t-1})=x_0x_1\cdots x_{t-1}$.
Then, it is easy to see that after doing pivot on any edge mentioned above, the obtained function  
does not fulfil the condition of definition~\ref{hyperpivot}.
\end{proof}


\begin{lem} 
Let $f\in {\cal F}^{n,t}$. Then the number of flat spectra of $f$ w.r.t. $\{I,H,N\}^n$ is at least 
$(n+1)(t+1)2^{n-t-1}$.
\end{lem}

\begin{proof} 
Let $f\in {\cal F}^{n,t}$. By theorem~\ref{family}, its number of flat spectra w.r.t. $\{I,H\}^n$ 
is at least $(t+1)2^{n-t-1}$; furthermore, we can see that all the flat spectra correspond to 
graph operations, so the resulting state is associated to a graph. It can be shown~\cite{RP:BC1} that 
local complementation at vertex $j$ is realised by the application of $N_j$ to the bipolar vector of 
the function, followed by a diagonal transform, which implies that the result of applying $N_j$ to 
the bipolar vector of a function associated to a (simple, undirected) graph is always flat (this also follows 
as a special case of theorem \ref{NAPFext}).  On the other 
hand, the result of applying the identity transform to the bipolar vector of a function associated to 
a graph is always flat. Therefore, the number of flat spectra of $f$ w.r.t. $\{I,H,N\}^n$ is at least 
$n+1$ times its number of flat spectra w.r.t. $\{I,H,N\}^n$; i.e. $(n+1)(t+1)2^{n-t-1}$.
\end{proof}


\section{Pivot Orbits and Codes}

A binary linear $[n,k]$ code $\mathcal{C}$ is a linear subspace of $\GF(2)^n$ of dimension
$k$. The $2^k$ elements of $\mathcal{C}$ are called \emph{codewords}.
We define the \emph{dual} of the code $\mathcal{C}$ with respect to the
standard inner product, $\mathcal{C}^\perp = 
\{ \boldsymbol{u} \in \GF(2)^n \mid \boldsymbol{u} \cdot \boldsymbol{c}=0, 
\forall \boldsymbol{c} \in \mathcal{C} \}$.
The code $\mathcal{C}$ can be defined by a $k \times n$ \emph{generator matrix}, $C$, whose 
rows span $\mathcal{C}$. 
Two codes, $\mathcal{C}$ and $\mathcal{C}'$, are considered to be \emph{equivalent} if
one can be obtained from the other by some permutation of the coordinates, or equivalently,
a permutation of the columns of the generator matrix.
$\mathcal{C}$~is called \emph{self-dual} if $\mathcal{C} = \mathcal{C}^\perp$,
and \emph{isodual} if $\mathcal{C}$ is equivalent to $\mathcal{C}^\perp$.
Self-dual and isodual codes must be \emph{even}, i.e., all codewords must have even weight.
A set of $k$ independent columns of $C$ is called an \emph{information set} of $\mathcal{C}$. 
The remaining $n-k$ columns is called a \emph{redundancy set}. 
We can permute the columns of $C$ such that an information set makes up the first $k$ columns. 
This matrix can now be transformed, by elementary row operations, into a matrix of the form 
$C' = (I \mid P)$, where $I$ is a $k \times k$
identity matrix, and $P$ is some $k \times (n-k)$ matrix. The matrix $C'$ generates a code
equivalent to $\mathcal{C}$ and is said to be of \emph{standard form}.
It follows that every code is equivalent to a code with generator matrix of standard form.
The matrix $H = (P^\text{T} \mid I)$, 
where $I$ is an $(n-k) \times (n-k)$ identity matrix is called the \emph{parity check matrix} of 
$\mathcal{C}$. Observe that $GH^\text{T} = \boldsymbol{0}$, where $\boldsymbol{0}$ is the 
all-zero vector. It follows that $H$ must be the generator matrix of $\mathcal{C}^\perp$.
A code is \emph{decomposable} if it can be written as the \emph{direct sum} of two smaller codes.
For example, let $\mathcal{C}$ be an $[n,k]$ code and $\mathcal{C}'$ an $[n',k']$ code. The
direct sum, $\mathcal{C} \oplus \mathcal{C}' = \{u||v \mid u \in \mathcal{C}, v \in \mathcal{C}'\}$,
where $||$ means concatenation, is an $[n+n',k+k']$ code.

It has previously been discovered that the LC orbits of simple undirected graphs corresponds to
the equivalence classes of \emph{self-dual additive codes over 
$\GF(4)$}~\cite{Bou:Graph,Glynn:Graph,nest2,selfdual}. We now show that pivot orbits of bipartite graphs correspond
the equivalence classes of binary linear codes.

\begin{df}\label{def:code}
Let $\mathcal{C}$ be a binary linear $[n,k]$ code. Let $C = (I \mid P)$ be a generator matrix
of standard form that generates a code equivalent to $\mathcal{C}$. Then the code
$\mathcal{C}$ corresponds to the $(k,n-k)$-bipartite graph on $n$ vertices with
adjacency matrix
\[
\Gamma = \begin{pmatrix}\boldsymbol{0}_{k\times k} & P \\ P^\text{T} & \boldsymbol{0}_{(n-k)\times (n-k)}\end{pmatrix},
\]
where $\boldsymbol{0}$ denote all-zero matrices of the specified dimensions.
Note that the graph corresponding to a code, like the generator matrix, is not uniquely defined. 
\end{df}

An alternative description of the relationship between bipartite graphs and
codes was given by Parker and Rijmen~\cite{PR:QE}.
We have previously shown how a graph corresponds to a Boolean function.
Applying the Hadamard transform, $H$, to all variables corresponding to vertices in one partition
of the graph (and $I$ to the other variables)  produces, to within normalisation, the binary
 \emph{indicator vector} 
of the corresponding code $\mathcal{C}$, i.e., a vector
$(s_{\boldsymbol{c}})$, $\boldsymbol{c} \in \mathbb{Z}_2^n$, where
$s_{\boldsymbol{c}} = 1$ if $\boldsymbol{c} \in \mathcal{C}$, and 
$s_{\boldsymbol{c}} = 0$ otherwise. More explicitly, for $s = (-1)^p$, and $p$ a quadratic Boolean function
representing the bipartite graph of the code $\mathcal{C}$, we have 
$(s_{\boldsymbol{c}}) = \mu(I \otimes \cdots \otimes I \otimes \cdots \otimes H \otimes \cdots \otimes H)s$, 
with $\mu$ some normalisation constant. Similarly, applying the $H$ transform to the
vertices of the other partition will give the indicator vector of $\mathcal{C}^\perp$

\begin{lem}\label{pivotP}
Let $G=(V,E)$ be a $(k,n-k)$-bipartite graph derived from the standard 
form generator matrix $C$ of the $[n,k]$ code $\mathcal{C}$. Let $G'$ be the graph obtained by 
pivoting on the edge $uv \in E$, followed by a swapping of vertices $u$ and $v$.
Both $G$ and $G'$ have generator matrices of the form given in definition~\ref{def:code},
with submatrices $P$ and $P'$, respectively.
Let the rows of $P$ be labelled $1,2,\ldots,k$, and let the columns of $P$
be labelled $k+1,k+2,\ldots,n$. 
Assuming, without loss of generality, that $u \le k$ and $v > k$,
$P$ can be transformed into $P'$ by the following steps.
\begin{enumerate}
\item Store the current value of column $v$ for later.
\item Add row $u$ to all rows in $\mathcal{N}(v) \setminus \{u\}$,
i.e., all rows that have 1 in coordinate $v$, except row $u$.
(Observe that column $v$ is now the basis vector $e_u$, i.e., it has 0 in all coordinates except coordinate $u$.)
\item Reset column $v$ to the value that was stored initially.
\end{enumerate}
\end{lem}
\begin{proof}
According to lemma~\ref{bip}, $G'$ will remain $(k,n-k)$-bipartite.
The transformation of $P'$ follows from definition~\ref{pivot2}.
Pivoting on the edge $ij$ of the bipartite graph $G$ is done
by `toggling' all pairs of vertices $\{x,y\}$, where $x \in \mathcal{N}(u) \setminus \{v\}$
and $y \in \mathcal{N}(v) \setminus \{u\}$.
This is obtained by step~2 above, since row $u$ of $P$ defines $\mathcal{N}(u)$,
and column $v$ defines $\mathcal{N}(v)$.
But in step~2 we have also `toggled' the pairs $\{u,y\}$, where $y \in \mathcal{N}(v) \setminus \{u\}$,
and we need steps~1 and 3 to correct this.
\end{proof}

\begin{theorem}\label{thm:pivotswap}
Let $G=(V,E)$ be the bipartite graph derived from the standard form generator matrix
$C=(I \mid P)$ of the code $\mathcal{C}$.
The graph $G'$ obtained by pivoting on the edge $uv \in E$ and
then swapping vertices $u$ and $v$
corresponds to the standard form generator matrix $C'=(I \mid P')$ of
the code $\mathcal{C}'$. The code $\mathcal{C}'$ is equivalent to $\mathcal{C}$,
and can be obtained by interchanging coordinates $u$ and $v$ of $\mathcal{C}$.
\end{theorem}
\begin{proof}
Assume that $u \le k$ and $v > k$. The effect of pivoting on the submatrix $P$ was described in lemma~\ref{pivotP}.
Now consider the following operations on $C=(I \mid P)$, where
rows are labelled $1,2,\ldots,k$, and columns are labelled $1,2,\ldots,n$.
\begin{enumerate}
\item Observe that column $u$ is the basis vector $e_u$.
\item Add row $u$ to all rows in $\mathcal{N}(v) \setminus \{u\}$.
\item Column $v$ is now the basis vector $e_u$, and column $u$ has the value that column $v$ had initially.
\item Swap columns $u$ and $v$.
\end{enumerate}
Comparing this with the algorithm for pivoting on $P$, it is easy to see that the resulting matrix is
$C'=(I \mid P')$, i.e., the generator matrix corresponding to $G'$.
The operations we have performed on $C$ preserve the equivalence of linear codes, namely
row additions and the swapping of columns $u$ and $v$.
\end{proof}

\begin{cor}
Applying any sequence of pivot operations to the graph $G$ corresponding to the code $\mathcal{C}$
will produce a graph corresponding to a code equivalent to $\mathcal{C}$.
\end{cor}

Consider a code $\mathcal{C}$. It is possible to go from a generator matrix 
of standard form, $C = (I \mid P)$, to
a generator matrix of standard form, $C'$, of any code equivalent to $\mathcal{C}$
by one of the $n!$ possible permutations of the columns of $C$.
More precisely, we can get from $C$ to $C'$ via a combination of the following operations.
\begin{enumerate}
\item Permuting the columns of $P$.
\item Permuting the columns of $I$, followed by the same permutation on the rows of $P$, 
      to restore standard form.
\item Swapping columns from $I$ with columns from $P$, such that the first $k$ columns of the
      generator matrix is an information set, followed by some row additions to restore standard form.
\end{enumerate}

\begin{theorem}\label{thm:alwayspivot}
Let $\mathcal{C}$ and $\mathcal{C}'$ be equivalent codes. Let $C$ and $C'$ be any matrices
of standard form generating $\mathcal{C}$ and $\mathcal{C}'$.
Let $G$ and $G'$ be the bipartite graphs corresponding to $C$ and $C'$.
$G'$ must be isomorphic to a graph obtained by performing some sequence of pivot operations on $G$.
\end{theorem}
\begin{proof}
$\mathcal{C}$ and $\mathcal{C}'$ must be related by a combination of the operations 1, 2, and 3 listed above.
It is easy to see that operations 1 and 2 applied to $G$ produce a graph isomorphic to $G$. It remains to 
prove that operation 3 always correspond to some sequence of pivot operations.
We know from theorem~\ref{thm:pivotswap} that swapping columns $u$ and $v$ of $C$, where
$u$ is part of $I$ and $v$ is part of $P$, corresponds to pivoting on the edge $uv$ of $G$ and then
swapping vertices $u$ and $v$.
When $uv$ is not an edge of $G$, we can not swap columns $u$ and $v$ of $C$ via pivoting.
In this case, coordinate $v$ of column $u$ is 0, and column $u$ is the basis vector $e_u$.
Swapping these columns would result in a generator matrix where the first $k$ columns have 0 
at coordinate $u$. These columns can not correspond to an information set.
It follows that if $uv$ is not an edge of $G$, swapping columns $u$ and $v$ is not a valid
operation of type 3 in the above list. Thus graph pivoting covers all possible permutations
that map standard form generator matrices of equivalent codes to each other.
\end{proof}

Let us now consider the labelled graphs in the pivot orbit of $G=(V,E)$, i.e.,
graph isomorphism is not considered when the pivot orbit is generated.
$G$ is the bipartite graph representing the code $\mathcal{C}$.
When we pivot on the edge $uv \in E$, without swapping vertices $u$ and $v$ afterwards,
the resulting adjacency matrix will not be of the type we saw in definition~\ref{def:code}.
We can think of $G$ as a graph corresponding to the information set
$\{1,2,\ldots,k\}$ of $\mathcal{C}$. Pivoting on the edge $uv \in E$, 
where $u \le k$ and $v > k$, produces a graph representing another information 
set of $\mathcal{C}$, namely $\{1,2,\ldots,k\} \setminus \{u\} \cup \{v\}$.
With this interpretation, the next corollary follows from theorem~\ref{thm:alwayspivot}.

\begin{cor}
Let $G$ be the bipartite graph representing the code $\mathcal{C}$.
Each labelled graph in the pivot orbit of $G$ corresponds to an information set of $\mathcal{C}$.
Moreover, the number of information sets of $\mathcal{C}$ equals the number of
labelled graphs in the pivot orbit of $G$.
\end{cor}

\section{Enumeration of Pivot Orbits}

We have previously classified all self-dual additive codes over $\GF(4)$ of length up 
to~12~\cite{DanThes,selfdual}, by classifying orbits of simple undirected graphs
with respect to local complementation and graph isomorphism.
In Table~\ref{tab:lcorbits}, the sequence $(i_n^{LC})$ gives the number of LC orbits 
of connected graphs on $n$ vertices, while $(t_n^{LC})$ gives the total number of 
LC orbits of graphs on $n$ vertices. A representative from each LC orbit is available 
at \url{http://www.ii.uib.no/~larsed/vncorbits/}.

\begin{table}
\centering
\caption{Numbers of LC Orbits of Graphs on $n$ Vertices}
\label{tab:lcorbits}
\begin{tabular}{ccccccccccccc}
\toprule
& \multicolumn{12}{c}{$n$}\\
\cmidrule(l){2-13}
 & 1 & 2 & 3 & 4 & 5 & 6 & 7 & 8 & 9 & 10 & 11 & 12 \\
\midrule
$i_n^{LC}$ & 1 & 1 & 1 & 2 &  4 & 11 & 26 & 101 & 440 & 3,132 & 40,457 & 1,274,068 \\
$t_n^{LC}$ & 1 & 2 & 3 & 6 & 11 & 26 & 59 & 182 & 675 & 3,990 & 45,144 & 1,323,363 \\
\bottomrule
\end{tabular}
\end{table}

By recursively applying pivot operations 
to all edges of a graph, whilst checking for graph isomorphism using the program 
\emph{nauty}~\cite{nauty}, we can quickly find all members of the pivot orbit. 
Let $\boldsymbol{G}_n$ be the set of all unlabelled simple undirected connected 
graphs on $n$ vertices. Let the set of all distinct pivot orbits of connected graphs 
on $n$ vertices is a partitioning of $\boldsymbol{G}_n$ into $i_n^P$ disjoint sets.
Our previous classification of the LC orbits of all graphs of up to 12 vertices helps 
us to classify pivot orbits, since it follows from definition~\ref{pivot}
that each LC orbit can be partitioned into some set of disjoint pivot orbits. 
We have used this fact to classify all pivot orbits of graphs on up to 12 vertices.
In Table~\ref{tab:pivorbits}, the sequence $(i_n^{P})$ gives the number of pivot orbits 
of connected graphs on $n$ vertices, while $(t_n^{P})$ gives the total number of 
pivot orbits of graphs on $n$ vertices.
A database containing one representative from each pivot orbit
can be found at \url{http://www.ii.uib.no/~larsed/pivot/}.

\begin{table}
\centering
\caption{Numbers of Pivot Orbits of Graphs on $n$ Vertices}
\label{tab:pivorbits}
\begin{tabular}{rrrrrrrrrrrrrrr}
\toprule
$n$ & $i_n^P$ & $t_n^P$ & $i_n^{P,B}$ & $t_n^{P,B}$ \\
\midrule
 1 &         1 &         1&      1 &      1  \\
 2 &         1 &         2&      1 &      2  \\
 3 &         2 &         4&      1 &      3  \\
 4 &         4 &         9&      2 &      6  \\
 5 &        10 &        21&      3 &     10  \\
 6 &        35 &        64&      8 &     22  \\
 7 &       134 &       218&     15 &     43  \\
 8 &       777 &     1,068&     43 &    104  \\
 9 &     6,702 &     8,038&    110 &    250  \\
10 &   104,825 &   114,188&    370 &    720  \\
11 & 3,370,317 & 3,493,965&  1,260 &  2,229  \\
12 & 231,557,290 & 235,176,097 &  5,366 &  8,361  \\
13 &           &          & 25,684 & 36,441  \\
\bottomrule
\end{tabular}
\end{table}

We are particularly interested in bipartite graphs, because of their connection to binary linear codes.
For the classification of the orbits of bipartite graphs with respect to pivot and graph
isomorphism, the following technique is helpful.
If $G$ is an $(a,b)$-bipartite graph, it has $2^a + 2^b -2$ possible \emph{extensions}.
Each extension is formed by adding a new vertex and joining it to all possible combinations 
of at least one of the old vertices.
Let $\boldsymbol{P}_n$ be a set containing one representative from each
pivot orbit of all connected bipartite graphs on $n$ vertices. The set $\boldsymbol{E}_{n}$ 
be formed by making all possible extensions of all graphs in $\boldsymbol{P}_{n-1}$.
It can then be shown that $\boldsymbol{P}_n \subset \boldsymbol{E}_{n}$, i.e., that
the set $\boldsymbol{E}_{n}$ will contain at least one representative from each
pivot orbit of connected bipartite graphs on $n$ vertices.
The set $\boldsymbol{E}_{n}$ will be much smaller than $\boldsymbol{G}_n$, so it will be more
efficient to search for a set of pivot orbit representatives within $\boldsymbol{E}_{n}$.

In Table~\ref{tab:pivorbits}, the sequence $(i_n^{P,B})$ gives the number of 
pivot orbits of connected bipartite graphs on $n$ vertices, and
$(t_n^{P,B})$ gives the total number of pivot orbits of bipartite graphs on $n$ vertices.
A database containing one representative from each of these orbits
can be found at \url{http://www.ii.uib.no/~larsed/pivot/}.

It follows from theorem~\ref{thm:alwayspivot} that the orbits of simple undirected graphs
with respect to pivot and graph isomorphism correspond to equivalence classes of
binary linear codes. Note that the codes $\mathcal{C}$ and $\mathcal{C}^\perp$ correspond 
to isomorphic graphs. This means that the pivot orbit of an $[n,k]$ code is simultaneously
the pivot orbit of a non-equivalent $[n,n-k]$ code, with the exception of isodual codes,
which are equivalent to their duals.

\begin{theorem}\label{prop:codenumbers}
Let $k \ne \frac{n}{2}$. Then the number of inequivalent binary linear $[n,k]$ codes, 
which is also the number of inequivalent $[n,n-k]$ codes, is equal to twice the number of
pivot orbits of $(n-k,k)$-bipartite graphs.

When $n$ is even and $k = \frac{n}{2}$, the number of inequivalent binary linear $[n,k]$ codes
is equal to twice the number of pivot orbits of $(k,k)$-bipartite graphs minus
the number of isodual codes of length $n$.
\end{theorem}

Note that if we only consider connected graphs on $n$ vertices,
we get the number of indecomposable codes of length $n$, $i_n^{C}$.
The total number of codes can easily be 
derived from the values of $(i_n^{C})$. Table~\ref{tab:binarycodes} gives the number of pivot 
orbits of connected bipartite graphs on $n$ vertices $(i_n^{P,B})$, the
number of indecomposable binary linear codes of length $n$ $(i_n^{C})$, and
the number of indecomposable isodual codes of length $n$ $(i_n^{C_{iso}})$.

\begin{table}[t]
\centering
\caption{Numbers of Pivot Orbits and Binary Linear Codes}
\label{tab:binarycodes}
\begin{tabular}{rrrrrrrrrrrrrrr}
\toprule
$n$ & $i_n^{P,B}$ & $i_n^{C}$ & $i_n^{C_{iso}}$ \\
\midrule
 1 &         1 &        1 &      - \\
 2 &         1 &        1 &      1 \\
 3 &         1 &        2 &      - \\
 4 &         2 &        3 &      1 \\
 5 &         3 &        6 &      - \\
 6 &         8 &       13 &      3 \\
 7 &        15 &       30 &      - \\
 8 &        43 &       76 &     10 \\
 9 &       110 &      220 &      - \\
10 &       370 &      700 &     40 \\
11 &     1,260 &    2,520 &      - \\
12 &     5,366 &   10,503 &    229 \\
13 &    25,684 &   51,368 &      - \\
\bottomrule
\end{tabular}
\end{table}

The number of linear codes of high length can be calculated by using
computer algebra tools~\cite{isometry}, and a complete classification has 
been carried out for codes of length up to~14~\cite{ostergard} by using
a different graph based approach. We hope, however, that our method will 
be more efficient than existing algorithms for classifying special types of codes.

Finally, we have also enumerated the orbits of labelled graphs
with respect to the pivot operation only, i.e., not considering graph isomorphism.
In Table~\ref{tab:pivlabel}, the sequence $(i_n^{P,L})$ gives the number of pivot orbits 
of connected labelled graphs on $n$ vertices, while $(t_n^{P,L})$ gives the total number of 
pivot orbits of labelled graphs on $n$ vertices. Similarly $(i_n^{P,B,L})$ and $(t_n^{P,B,L})$
give the numbers for connected and unconnected bipartite labelled graphs.

\begin{table}[t]
\centering
\caption{Numbers of Pivot Orbits of Labelled Graphs on $n$ Vertices}
\label{tab:pivlabel}
\begin{tabular}{rrrrrrrrrrrrrrr}
\toprule
$n$ & $i_n^{P,L}$ & $t_n^{P,L}$ & $i_n^{P,B,L}$ & $t_n^{P,B,L}$ \\
\midrule
 1 &         1 &        1 &      1 &      1  \\
 2 &         1 &        2 &      1 &      2  \\
 3 &         2 &        6 &      1 &      5  \\
 4 &        11 &       29 &      4 &     18  \\
 5 &       119 &      240 &     26 &     92  \\
 6 &     2,303 &    3,623 &    251 &    693  \\
 7 &    80,923 &  105,564 &  3,412 &  7,613  \\
\bottomrule
\end{tabular}
\end{table}

\pagebreak

\end{document}